\documentclass[12pt]{article}
\usepackage{amsmath}
\usepackage{graphicx}
\usepackage{hyperref}
\usepackage[latin1]{inputenc}
\usepackage{mathtools}
\usepackage{amsfonts}
\usepackage{amssymb}
\usepackage[T1]{fontenc}
\usepackage{amsthm}

\theoremstyle{definition}
\newtheorem{definition}{Definition}[section]
 
\theoremstyle{remark}

\theoremstyle{theorem}
\newtheorem{theorem}{Theorem}[section]
\newtheorem{corollary}{Corollary}[theorem]
\newtheorem{lemma}[theorem]{Lemma}
\newtheorem*{remark}{Remark} 

\title{Sufficient condition for existence of special type of primitive normal elements over finite fields}
\author{Himangshu Hazarika and Dhiren Kumar Basnet*\\Department of Mathematical Sciences\\ Tezpur University, Napaam}
\date{}

\begin{document}
\maketitle

  \par
 \textbf{Abstract:} Let $\mathbb{F}_{q^n}$ be the extension of the field $\mathbb{F}_q$ of degree n, where $q$ is power of prime $p$, i.e $q=p^k$, where k is a positive integer. In this paper, we provide sufficient condition for the existence of a primitive normal element $\alpha\in\mathbb{F}_{q^n} $ such that $\alpha^2+\alpha+1$ is also primitive normal element over $\mathbb{F}_{q^n}$.   
 
\vspace{.5cm}

\par
\textbf{Keywords:} Finite field, primitive element, normal element, Character.
 \\
 \\
\section{Introduction}
\par Let $\mathbb{F}_q$ be a finite field of order $q=p^k$, where $p$ is some prime and $k$ is some positive integer. $\mathbb F_{q^n}$ be the extension field of $\mathbb{F}_{q}$ of degree $n$. For any finite field its multiplicative group $\mathbb F^*_{q^n}$ is cyclic.The generators of $\mathbb F^*_{q^n}$ are called $\mathit{primitive \thinspace elements}$ of $\mathbb F_{q^n}$. Any field of order $q^n$ (i.e.,  $\mathbb{F}_{q^n}$) has $\phi(q^n-1)$ primitive elements, where $\phi$ is the Euler's phi-function. An element $\alpha\in \mathbb{F}_{q^n}$ is called a $\mathit{normal \thinspace element}$ of $\mathbb{F}_{q^n}$ over $\mathbb{F}_{q}$ if $\{ \alpha,\alpha^q,\dots,\alpha^{q^{n-1}}\}$ is a basis of $\mathbb{F}_{q^n}(\mathbb{F}_{q})$. This basis is called $\mathit{normal\thinspace basis}$. Normal bases are quite useful in finite field theory, as they are applicable in coding theory, cryptography etc \cite{1}.It is well known (\cite{13},Theorem 2.35), that there exists a normal basis of $\mathbb{F}_{q^n}$ over $\mathbb{F}_q$. An element $\alpha\in\mathbb F^*_{q^n}$ is called a $\mathit{primitive \thinspace normal \thinspace element}$ if it is primitive as well as normal element of $\mathbb{F}_{q^n}$ over $\mathbb{F}_{q}$.   
\par
  Using the properties of primitive elements , modern day cryptosystems such as ElGamal crypto-\\system,$\mathit{ The \thinspace Diffie-}\,\mathit{Hellman\thinspace key\thinspace agreement\thinspace protocol}$, RSA cryptosystem are developed. Interestingly, even though there are $\phi(q^n-1)$ primitive elements in a finite field $\mathbb{F}_{q^n}$, but finding one such primitive element may be difficult, as there is no polynomial time algorithm to compute a primitive element. But one can determine a primitive element in context of other. The main focus is to prove the existence of primitive element in terms of another primitive element, thus making a choice between them for further applications. Same is applicable for primitive normal element. In our paper we establish a sufficient condition for existence of a primitive normal element in context of another one. Lenstra and Schoof\cite{11} proved the existence of primitive normal element for any finite field $\mathbb F_{q^n}$ over $\mathbb F_{q}$. Later, Cohen and Huczynska\cite{4} gave a computer free proof of the existence of a primitive normal element.
 \par 
 For any primitive element $\alpha \in \mathbb{F}_{q^n}$ and for a rational function $f$, $f(\alpha)$ need not be primitive in $\mathbb{F}_{q^n}$, for example, if we consider $f(x)=x+1$ over the field $\mathbb{F}_2$, then we see that $f(1)$ is not a primitive element in $\alpha \in \mathbb{F}_2$ whereas $1$ is a primitive element of $\mathbb{F}_2$. But for $f(x)=\frac{1}{x}$, $f(\alpha)$ is always primitive for $\alpha\in\mathbb{F}_{q^n}$. Cohen and Han studied the existence of primitive element $\alpha$ such that $f(\alpha)=\alpha+\alpha^{-1}$ is also primitive in a finite field. In 2012, Wang\cite{16} established a sufficient condition for existence of a primitive element $\alpha$ such that $\alpha+\alpha^{-1}$ is also primitive for the case $2|q$, later on it generalised by Leao\cite{17}. Further Tian and Qi\cite{18} proved the existence of a normal element $\alpha\in \mathbb{F}_{q^n}$ such that $\alpha^{-1}$ is also normal in $\mathbb{F}_{q^n}$ over $\mathbb{F}_{q}$ for $n\geq32$. Later Cohen and Huczynska\cite{5} proved the existence of a primitive normal element $\alpha\in \mathbb{F}_{q^n}$ such that $\alpha^{-1}$ is also primitive normal over $\mathbb{F}_{q}$ for $n\geq2$, except when $(q,n)$ is one of the pairs (2,3),  (2,4), (3,4), (4,3), (5,4). In 2018, Anju and R.K.Sharma\cite{14} established a sufficient condition for the existence of a primitive element $\alpha\in \mathbb{F}_{q^n}$, such that for $\alpha^2+\alpha+1$ is also primitive in $\mathbb{F}_{q^n}$. Further they established a sufficient condition for the existence of a primitive normal element $\alpha\in \mathbb{F}_{q^n}$, such that $\alpha^2+\alpha+1$ is also primitive in $\mathbb{F}_{q^n}$. 
 \par 
 In this paper we use the notation $\mathfrak{M}$ for the set $(q,n)$ such that $\mathbb{F}_{q^n}$ contains primitive normal element $\alpha$ such that $\alpha^2+\alpha+1$ is also primitive normal element in $\mathbb{F}_{q^n}$. For any positive integer $m>1$ and any $g(x)\in \mathbb{F}_q[x]$, 
 $\omega(m)$ and $\Omega_q(g)$ denote the number of prime divisors of $m$ and the number of monic irreducible divisors of $g$ over $\mathbb{F}_{q}$ respectively.

\vspace{.4cm} 
\section{Prerequisites}
First of all, we recall some definitions.
 \begin{definition} \textbf{Character}
 Let $G$ be a finite abelian group and $S := \{z\in \mathbb{C}: |z|=1\}$ be the multiplicative group of all complex numbers with modulus 1 . Then a character $\chi$ of $G$ is a homomorphism from $G$ into the group $S$, i.e $\chi(a_1a_2)=\chi(a_1)\chi(a_2)$ for all $a_1,a_2\in G$. The characters of $G$ forms a group under multiplication called $\mathit{dual\thinspace group}$ or $\mathit{character\thinspace group}$ of $G$ which is denoted by $\widehat{G}$. It is well known that $\widehat{G}$ is isomorphic to $G$. Again the character $\chi_0$ is denoted for the trivial character of $G$ defined as $\chi_0(a)=1$ for all $a \in G$.
 \par In a finite field $\mathbb{F}_{q^n}$, there are two types of abelian groups, namely additive group $\mathbb{F}_{q^n}$ and multiplicative group $\mathbb{F}^*_{q^n}$. So, there are two types of characters of a finite field $\mathbb{F}_{q^n}$, namely $\mathit{additive \thinspace character }$ for $\mathbb{F}_{q^n}$ and $\mathit{multiplicative \thinspace character }$ for $\mathbb{F}^*_{q^n}$. Multiplicative characters are extended from $\mathbb{F}^*_{q^n}$ to $\mathbb{F}_{q^n}$ by the  rule
 \hspace{.1cm}  $\chi(0)=\begin{cases}
                         0 \,\mbox{if}\, \chi\neq\chi_0\\
                         1 \,\mbox{if}\, \chi=\chi_0 
                         \end{cases} $ 
                         
  \par Since $\widehat{\mathbb{F}^*_{q^n}} \cong \mathbb{F}^*_{q^n}$, so $\widehat{\mathbb{F}^*_{q^n}}$ is cyclic and for any divisor $d$ of $q^n-1$, there are exactly $\phi(d)$ characters of order $d$ in $\widehat{\mathbb{F}^*_{q^n}}$.
  \end{definition}
  \begin{definition} \textbf{e-free element}
  For any divisor $e$ of $q^n-1$, an element $\alpha\in \mathbb{F}_{q^n}$ is called $\mathit{e-free}$, if for any $d|e, \alpha=\beta^d$ where $\beta\in \mathbb{F}_{q^n}$ implies $d=1 $ i.e, if gcd$(d,\frac{q^n-1}{ord_{q^n}(\alpha)})=1$. Hence an element $\alpha\in \mathbb{F}^*_{q^n}$ is primitive if and only if it is $q^n-1$-free. 
  \end{definition} 
\begin{definition} \textbf{$\mathbb{F}_q$-order of an element}
 \\ The additive group of $\mathbb{F}_{q^n}$ is a $\mathbb{F}_{q}[x]$-module under the rule
 \par \hspace{3cm}$fo\alpha=\overset{m}{\underset{i=1}{\sum}} a_i\alpha^{q^i}$; for $\alpha\in \mathbb{F}_{q^n}$ and 
 $f(x)= \overset{m}{\underset{i=1}{\sum}}a_ix^i\thinspace \in \mathbb{F}_{q}[x]$.
 \\ For $\alpha \in \mathbb{F}_{q^n}$, the $\mathbb{F}_q$-order of $\alpha$ is the monic $\mathbb{F}_q$-divisor $g$ of $x^n-1$ of minimal degree such that $g o\alpha=0$.
 \end{definition}
 \begin{definition} \textbf{$g$-free element}
 \\ Let $g$ be a divisor of $x^n-1$. If, $\alpha=ho\beta$ where $\beta\in \mathbb{F}_{q^n}$, $h$ is a divisor of $x^n-1$ implies $h=1$, then $\alpha$ is called $g$-free in $\mathbb{F}_{q^n}$. Hence an element $\alpha\in \mathbb{F}_{q^n}$ is normal if and only if it is $x^n-1$ free. 
 \end{definition} 
\begin{definition}\textbf{Character function} 
\\For any $e|q^n-1$, Cohen and Huczynska\cite{4,5} defined  the character function for the subset of $e$-free elements of $\mathbb{F}^*_{q^n}$ by
\par \hspace{4cm }$\rho_e: \alpha\mapsto\theta(e)\underset{d|e}{\sum}(\frac{\mu(d)}{\phi(d)}\underset{\chi_d}{\sum}\chi_d(\alpha))$
    \\ where $\theta(e):=\frac{\phi(e)}{e}$, $\mu$ is the M\"obius function and $\chi_d$ stands for any multiplicative character of order $d$.
    \par Again, for any monic $\mathbb{F}_q$-divisor $g$ of $x^n-1$, a typical additive character $\psi_g$ of $\mathbb{F}_q$-order $g$ is one such that $\psi_g o g$ is the trivial character in $\mathbb{F}_{q^n}$ and $g$ is of minimal degree satisfying this property. 
    \\Furthermore, there are $\Phi_q(g)$ characters $\psi_g$, where $\Phi_q(g)= (\mathbb{F}_q[x]/g\mathbb{F}_q[x])^*$ is the analogue of Euler function on $\mathbb{F}_{q}[x]$. 
    \\ Then the character function for the set of $g$-free elements in $\mathbb{F}_{q^n}$, for any $g|x^n-1$ is given by
    \par \hspace{4cm }$\kappa_g :\alpha\mapsto\Theta(g)\underset{f|g}{\sum}(\frac{\mu^\prime(f)}{\Phi(f)}\underset{\psi_f}{\sum}\psi_f(\alpha))$
    \\when $\Theta(g):= \frac{\Theta(g)}{q^{deg(g)}}$, sum runs over all additive characters $\psi_f$ of $\mathbb{F}_q$-order g and $\mu^\prime$ is the analogue of the M\"obius function which is defined as follows: \\
     \hspace{.1cm}  $\mu^\prime(g)=\begin{cases}
                         (-1)^s\hspace{.2cm} \mbox{if \,g \, is \, a \, product\, of \, s\, distinct\, monic\, polynomials} \\
                        \hspace{.1cm} 0 \qquad\mbox{otherwise}\\ 
                         \end{cases} $ 
    
From Cohen and Huczynska\cite{5}, we have the following about the typical additive character.
\\ Let $\lambda$ be the canonical additive character of $\mathbb{F}_q$. Thus for $\alpha\in \mathbb{F}_q$ this character is defined as 
\par \hspace{5cm} $\lambda(\alpha)= \exp^{2\pi iTr(\alpha)/p}$   
\\ where $Tr(\alpha)$ is absolute trace of $\alpha$ over $\mathbb{F}_p$.
\\ Now let $\psi_0$ be canonical additive character of $\mathbb{F}_{q^n}$, it is simply the lift of $\lambda$ to $\mathbb{F}_{q^n}$ i.e., $\psi_0(\alpha)=\lambda(Tr(\alpha)), \, \alpha\in 
\mathbb{F}_{q^n}$. Now for any $\delta\in \mathbb{F}_{q^n}$, let $\psi_\delta$ be the character defined by $\psi_\delta(\alpha)=\psi_0(\delta\alpha), \, \alpha\in \mathbb{F}_{q^n}$.
Define the subset $\Delta_g$ of $\mathbb{F}_{q^n}$ as the set of $\delta$ for which $\psi_\delta$ has $\mathbb{F}_{q}$-order $g$. So we may also write $\psi_{\delta_g}$ for $\psi_\delta$, where $\delta_g\in \Delta_g$. So with the help of this we can express any typical additive character $\psi_g$ in terms of $\psi_{\delta_g}$ and further we can express this in terms of canonical additive character $\psi_0$.
 
\end{definition}
\section{Theorems and lemmas used in this paper}
In this section we recall some theorems which will be used throughout our discussions 
\begin{theorem} ${\mathbf{(\cite{13}, Theorem\, 5.4)}}$
If $\chi$ is any nontrivial character of a finite abelian group $G$ and $\alpha\in G$ any nontrivial element then 
\par \hspace{1cm}  $\underset{\alpha\in G}{\sum}\chi(\alpha)=0$\hspace{.5cm} and \hspace{.5cm}$\underset{\chi\in \widehat{G}}{\sum}\chi(\alpha)=0$. 
\end{theorem}
\begin{theorem}
$\mathbf{(\cite{13}, theorem\, 5.11)}$ Let $\chi$ be a nontrivial multiplicative character and $\psi$ a nontrivial additive character of $\mathbb{F}_{q^n}$. Then
\par \hspace{5cm}$|\underset{\alpha\in \mathbb{F}^*_{q^n}}{\sum}\chi(\alpha)\psi(\alpha) |=q^{n/2}$.
\end{theorem}
\begin{theorem} $\mathbf{(\cite{7}, Corollary\, 2.3.})$
\\ Consider any two nontrivial multiplicative characters $\chi_1,\chi_2$ of the finite field $\mathbb{F}_{q^n}$. Again, let $f_1(x)$ and $f_2(x)$ be two monic pairwise co-prime polynomials in $\mathbb{F}_{q^n}[x]$, such that none of $f_i(x)$ is of the form $g(x)^{ord(\chi_i)}$ for $i=1,2$; where $g(x)\in \mathbb{F}_{q^n}[x]$ with degree at least 1. Then 
\par \hspace{3cm} $|\underset{\alpha\in \mathbb{F}_{q^n}}{\sum}\chi_1(f_1(\alpha))\chi_2(f_2(\alpha))|\leq (n_1+n_2-1)q^{n/2}$
\\ where $n_1$ and $n_2$ are the degrees of largest square free divisors of $f_1$ and $f_2$ respectively.
\end{theorem}
\begin{theorem}
$\mathbf{Weil's\, Theorem\,(\cite{13}, Theorem\, 5.38)}$ 
\\ Let $f\in \mathbb{F}_{q}[x]$ of degree $n\geq1$ with gcd$(n,q)=1$ and $\psi$ be nontrivial additive character of $\mathbb{F}_{q^n}$. Then
\par \hspace{4cm} $\vert\underset{\alpha\in\mathbb{F}_{q^n}}{\sum}\psi(f(\alpha))|\leq \,(n-1)q^{n/2}$.
\end{theorem}
\begin{theorem}$\mathbf{(\cite{13}, Theorem\, 5.41)}$
\\Let $\chi$ be a multiplicative character of $\mathbb{F}_{q^n}$ of order $m>1$ and $f\in \mathbb{F}_{q^n}[x]$ be a monic polynomial of positive degree that is not an $m^{th}$ power of a polynomial over $\mathbb{F}_{q^n}$. Let $d$ be the number of distinct roots of $f$ in its splitting field over $\mathbb{F}_{q^n}$. Then for every $a\in \mathbb{F}_{q^n}$, we have
\par \hspace{4cm} $ |\underset{\alpha\in\mathbb{F}_{q^n}}{\sum}\chi(af(\alpha))|\leq(d-1)q^{n/2}$
\end{theorem}
\begin{theorem} $\mathbf{(\cite{2}, Theorem\, 5.6)}$
\\ Let $f_1(x), f_2(x),\ldots, f_k(x)\in \mathbb{F}_{q^n}[x]$ be distinct irreducible polynomials and $g(x)$ be rational function over $\mathbb{F}_{q^n}$. Let $\chi_1,\chi_2.\ldots,\chi_k$ be multiplicative characters and $\psi$ be a nontrivial additive character of $\mathbb{F}_{q^n}$. Suppose that $g(x)$ is not of the form $r(x)^q-r(x)$ in $\mathbb{F}_{q^n}[x]$. Then
\par $\left | \underset{\underset{f_i(\alpha)\neq0, g(\alpha)\neq\infty}{\alpha\in\mathbb{F}_{q^n}}}{\sum} \chi_1(f_1(\alpha))\chi_2(f_2(\alpha))\ldots\chi_k(f_k(\alpha))\psi(g(\alpha))\right|$
\par\hspace{7cm}$\leq (n_1+n_2+n_3+n_4-1)q^{n/2}$ 
\\where $n_1= \overset{k}{\underset{j=1}{\sum}}deg(f_j),\, n_2= \mbox{max}(\mbox{deg}(g),0),\, n_3$ is the degree of denominator of $g(x)$ and $n_4$ is sum of degrees of those irreducible polynomials dividing the denominator of $g$, but distinct from $f_j(x),\, j=1,2,\dots,k$. 
\end{theorem}
\begin{lemma} $\mathbf{(\cite{11}, Lemma\, 2.6)}$
\\ Let $n>1, l>1$ be integers and $\Lambda$ be the set of primes $\leq\,l$. Set $L:= \underset{r\in\Lambda}{\Pi}r$. Assume that every prime factor $r<l$ of $n$ is contained in $\Lambda$. Then
\par \hspace{4cm} $\omega(n)\leq\frac{log\, n-log\,L}{log \,l}+|\Lambda|$ \hfill (3.1)
\\ Let $m$ be a positive integer and $p_m$  be the $m^{th}$ prime. Now we can take $l=p_m$, and then $\Lambda$ is the set of primes no more than $p_m$ , $|\Lambda|=m$ i.e., so the inequality (3.1)
becomes 
\par \hspace{4cm}$\omega(n)\leq\frac{log\, n-\overset{m}{\underset{i=1}{\sum}}\log\,p_i}{log \,p_m}+ \,m$ \hfill (3.2)
\end{lemma}
\begin{lemma} $\mathbf{(\cite{14}, Lemma \, 2.7)}$
\\ Let $q$ be a prime power and $n$ be a positive integer. Let $\Omega := \Omega_q(x^n-1)$. Then we have \\ $\Omega\leq \{n\,+\, \mbox{gcd}(n,q-1)\}/2$. In particular, $\Omega\leq n$ and $\Omega=n$ if and only if $n|q-1$. Moreover, $\Omega\leq\frac{3}{4}n$ if $n\nmid q-1$.
\end{lemma}
\begin{lemma} $\mathbf{(\cite{14}, Lemma \, 3.1)}$
\\ For any positive integer $N,\, 2^{\omega(N)}< \, C(N) \,N^{1/5}$, where $C(N)<11.25$.
\par \hspace{2cm} Moreover 
\hspace{.1cm}  $ C(N)<\begin{cases}
                         7.77\; \mbox{if}\; 5\nmid N\\
                         8.31\; \mbox{if}\; 7\nmid N\\ 
                         \end{cases} $

\end{lemma}
\vspace{1.5cm}
\section{ Main results}
 Let $N_{q^n}(m_1,m_2,g_1,g_2)$ be the number of $\alpha\in \mathbb{F}_{q^n}$, such that $\alpha$ is $m_1$-free, $\alpha^2+\alpha+1$ is $m_2$-free, $\alpha$ is $g_1$-free and $\alpha^2+\alpha+1$ is $g_2$-free, where $m_1,m_2$ are positive integers and $g_1, g_2$ are any polynomials over $\mathbb{F}_q$. We use the notations $\chi_1$ and $\psi_1$ to denote the trivial multiplicative and additive characters respectively.

\vspace{1cm}
\begin{theorem}

 Let $q=p^k$ for some prime $p\neq 2,3$; $k\in\mathbb{N}$ and  $n$ be a positive integer. Let us write $\omega := \omega(q^n-1)$ and $\Omega := \Omega_q(x^n-1)$.  If $q^{n/2}> 4.2^{2\omega+2\Omega}$, then $(q,n)\in \mathfrak{M}$.
\end{theorem}

\vspace{.3cm}
\textbf{Proof:} By definition 
 \begin{align*}
 N_{q^n}(q^n-1,q^n-1,x^n-1,x^n-1)&= \underset{\alpha\in\mathbb{F}^*_{q^n}}{\sum}\rho_{q^n-1}(\alpha)\rho_{q^n-1}(\alpha^2+\alpha+1)\kappa_{x^n-1}(\alpha)\kappa_{x^n-1}(\alpha^2+\alpha+1)\\
 &= \theta(q^n-1)^2\Theta(x^n-1)^2\underset{\alpha\in \mathbb{F}^*_{q^n}}{\sum}\underset{d,h|q^n-1}{\sum}\,\underset{g,f|x^n-1}{\sum}\frac{\mu(d)\mu(h)\mu^\prime(g)\mu^\prime(f)}{\phi(d)\phi(h)\Phi(g)\Phi(f)}\\
 & \quad \underset{\chi_d,\chi_h}{\sum}\underset{\psi_g,\psi_f}{\sum}\chi_d(\alpha)\chi_h(\alpha^2+\alpha+1)\psi_g(\alpha)\psi_f(\alpha^2+\alpha+1)\\
 &= \theta(q^n-1)^2\Theta(x^n-1)^2(\overset{16}{\underset{i=1}{\sum}}S_i) \\
\end{align*}
\vspace{.5cm}
If  $S_1$ is taken over $d=h=1=g=f$,then 
\begin{align*}
S_1 &= \underset{\alpha\in \mathbb{F}^*_{q^n}}{\sum}\underset{d=1=h}{\sum}\,\underset{g=1=f}{\sum}\frac{\mu(d)\mu(h)\mu^\prime(g)\mu^\prime(f)}{\phi(d)\phi(h)\Phi(g)\Phi(f)} \underset{\chi_d,\chi_h}{\sum}\underset{\psi_g,\psi_f}{\sum}\chi_d(\alpha)\chi_h(\alpha^2+\alpha+1)\psi_g(\alpha)\psi_f(\alpha^2+\alpha+1)\\
&= \underset{\alpha\in \mathbb{F}^*_{q^n}}{\sum}\underset{d=1=h}{\sum}\,\underset{g=1=f}{\sum}\left (\frac{\mu(1)\mu^\prime(1)}{\phi(1)\Phi(1)}\right )^2 \underset{\chi_1,\chi_1}{\sum}\underset{\psi_1,\psi_1}{\sum}\chi_1(\alpha)\chi_1(\alpha^2+\alpha+1)\psi_1(\alpha)\psi_1(\alpha^2+\alpha+1)\\
&= \underset{\alpha\in\mathbb{F}^*_{q^n}}{\sum}1 = q^n-1\\
\end{align*}
If $S_2$ is taken over $d\neq1, h=1=g=f$, then\\
\begin{align*}
|S_2| &= \left|\underset{\alpha\in \mathbb{F}^*_{q^n}}{\sum}\underset{1\neq d|q^n-1}{\sum}\,\underset{g=1=f}{\sum}\frac{\mu(d)}{\phi(d)} \underset{\chi_d}{\sum}\underset{\psi_1,\psi_1}{\sum}\chi_d(\alpha)\chi_1(\alpha^2+\alpha+1)\psi_1(\alpha)\psi_1(\alpha^2+\alpha+1) \right | \\
&= \left | \underset{\alpha\in \mathbb{F}^*_{q^n}}{\sum}\underset{1\neq d|q^n-1}{\sum}\frac{\mu(d)}{\phi(d)}\underset{\chi_d}{\sum}\chi_d(\alpha)\right | \\
\end{align*}
\par \hspace{1.5cm}$\leq \underset{1\neq d|q^n-1}{\sum}\frac{\mu(d)}{\phi(d)}\underset{\chi_d}{\sum}\left |\underset{\alpha\in \mathbb{F}^*_{q^n}}{\sum}\chi_d(\alpha)\right |$ \\

 By theorem 3.1, we have $|S_2|= 0$\\
\\If $S_3$ is taken over $h\neq 1, d=1=g=f$, then\\
\begin{align*}
|S_3|&= \left|\underset{\alpha\in \mathbb{F}^*_{q^n}}{\sum}\underset{1\neq h|q^n-1}{\sum}\,\underset{g=1=f}{\sum}\frac{\mu(h)}{\phi(h)} \underset{\chi_h}{\sum}\underset{\psi_1,\psi_1}{\sum}\chi_1(\alpha)\chi_h(\alpha^2+\alpha+1)\psi_1(\alpha)\psi_1(\alpha^2+\alpha+1) \right |\\
&\leq \underset{\underset{h\, squarefree}{1\neq h|q^n-1}}{\sum}\frac{1}{\phi(h)}\underset{\chi_h}{\sum}\left| \underset{\alpha\in \mathbb{F}_{q^n}}{\sum}\chi_h(\alpha^2+\alpha+1)-\chi_h(1)\right|
\end{align*}
By theorem 3.5, we have  $\left |\underset{\alpha\in\mathbb{F}_{q^n}}{\sum}\chi_h(\alpha^2+\alpha+1)\right|\leq q^{n/2}.\; \mbox{Using this, } \underset{\chi_h}{\sum}1=\phi(h) \; \mbox{and}\;\\ \underset{\underset{h\, squarefree}{1\neq h|q^n-1}}{\sum}1= 2^\omega-1, \; \mbox{we get} \, \, |S_3|\leq (q^{n/2}+1)(2^\omega-1)$
\vspace{.3cm}
\\If $S_4$ is taken over $d\neq1,\, h\neq1,\, g=1=f$, then 
\begin{align*}
|S_4| &= \left|\underset{\alpha\in \mathbb{F}^*_{q^n}}{\sum}\underset{1\neq d,h|q^n-1}{\sum}\,\underset{g=1=f}{\sum}\frac{\mu(d)\mu(h)}{\phi(d)\phi(h)} \underset{\chi_d,\chi_h}{\sum}\underset{\psi_1,\psi_1}{\sum}\chi_d(\alpha)\chi_h(\alpha^2+\alpha+1)\psi_1(\alpha)\psi_1(\alpha^2+\alpha+1) \right | \\
&\leq \underset{\underset{d,h \, squarefree}{1\neq d,h|q^n-1}}{\sum}\frac{1}{\phi(d)\phi(h)}\underset{\chi_d,\chi_h}{\sum}\left | \underset{\alpha\in\mathbb{F}^*_{q^n}}{\sum}\chi_d(\alpha)\chi_h(\alpha^2+\alpha+1)\right|
\end{align*}
 By theorem 3.3, we have\\ 

$|S_4|\leq \underset{\underset{d,h \, squarefree}{1\neq d,h|q^n-1}}{\sum}\frac{1}{\phi(d)\phi(h)}\underset{\chi_d,\chi_h}{\sum}2q^{n/2}\;= 2q^{n/2}(2^\omega-1)^2$\\
\\If $S_5$ is taken over $d=1=h, g\neq1, f=1$, then
\begin{align*}
|S_5| &=\left | \underset{\alpha\in\mathbb{F}^*_{q^n}}{\sum}\underset{1\neq g|x^n-1}{\sum}\frac{\mu^\prime(g)}{\Phi(g)}\underset{\psi_g}{\sum}\psi_g(\alpha)\right|\\
&\leq \underset{\underset{g\, squarefree}{1\neq g|x^n-1}}{\sum}\frac{1}{\Phi(g)}\underset{\psi_g}{\sum}\left| \underset{\alpha\in \mathbb{F}_{q^n}}{\sum}\psi_g(\alpha)-\psi_g(0)\right|
\end{align*}
 Now by applying theorem 3.1 and $\psi_g(0)=1$,  we have
$|S_5|\leq \underset{\underset{g\, squarefree}{1\neq g|x^n-1}}{\sum}\frac{1}{\Phi(g)}\underset{\psi_g}{\sum}1$\\
 Then by using the facts $ \underset{\psi_g}{\sum}1=\Psi(g) \, \,
\mbox{and}  \underset{\underset{g\, squarefree}{1\neq g|x^n-1}}{\sum}1=2^\Omega-1 \, \, 
\mbox{we have}\; |S_5|\leq (2^\Omega-1)$\\
If $S_6$ is taken over $h=1,d\neq1,g\neq1,f=1$, then
\begin{align*}
|S_6|&=\left | \underset{\alpha\in\mathbb{F}^*_{q^n}}{\sum}\underset{1\neq d|q^n-1}{\sum}\underset{1\neq g|x^n-1}{\sum}\frac{\mu(d)\mu^\prime(g)}{\phi(d)\Phi(g)}\underset{\chi_d}{\sum}\underset{\psi_g}{\sum}\chi_d(\alpha)\psi_g(\alpha)\right |\\
&\leq \underset{\underset{d\, square\,free}{1\neq d|q^n-1}}{\sum}\;\underset{\underset{g\, square\,free}{1\neq g|x^n-1}}{\sum}\frac{1}{\phi(d)\Phi(g)}\underset{\chi_d}{\sum}\underset{\psi_g}{\sum}\left |\underset{\alpha\in\mathbb{F}^*_{q^n}}{\sum}\chi_d(\alpha)\psi_g(\alpha)\right |
\end{align*}
Using theorem 3.2, we have\, $\left |\underset{\alpha\in\mathbb{F}^*_{q^n}}{\sum}\chi_d(\alpha)\psi_g(\alpha)\right |\leq q^{n/2}$,\, 
and hence $ |S_6|\leq q^{n/2}(2^\omega-1)(2^\Omega-1)$\\

If $S_7$ is taken over $d=1,h\neq1,g\neq1, f=1$,then

\begin{align*}
|S_7|&=\left | \underset{\alpha\in\mathbb{F}^*_{q^n}}{\sum}\underset{1\neq h|q^n-1}{\sum}\underset{1\neq g|x^n-1}{\sum}\frac{\mu(h)\mu^\prime(g)}{\phi(h)\Phi(g)}\underset{\chi_h}{\sum}\underset{\psi_g}{\sum}\chi_h(\alpha^2+\alpha+1)\psi_g(\alpha)\right |\\
&\leq \underset{\underset{h\, square\,free}{1\neq h|q^n-1}}{\sum}\;\underset{\underset{g\, square\,free}{1\neq g|x^n-1}}{\sum}\frac{1}{\phi(h)\Phi(g)}\underset{\chi_h}{\sum}\underset{\psi_g}{\sum}\left |\underset{\alpha\in\mathbb{F}^*_{q^n}}{\sum}\chi_h(\alpha^2+\alpha+1)\psi_g(\alpha)\right |\\
&\leq \underset{\underset{h\, square\,free}{1\neq h|q^n-1}}{\sum}\;\underset{\underset{g\, square\,free}{1\neq g|x^n-1}}{\sum}\frac{1}{\phi(h)\Phi(g)}\underset{\chi_h}{\sum}\underset{\psi_g}{\sum}\left |\underset{\alpha\in\mathbb{F}_{q^n}}{\sum}\chi_h(\alpha^2+\alpha+1)\psi_g(\alpha)-\chi_h(1)\psi_g(0)\right |
\end{align*}
By applying theorem 3.6, we have$\left |\underset{\alpha\in\mathbb{F}^*_{q^n}}{\sum}\chi_h(\alpha^2+\alpha+1)\psi_g(\alpha)\right |\leq(2q^{n/2}+1)$ and hence \\
$ |S_7|\leq(2q^{n/2}+1)(2^\omega-1)(2^\Omega-1)$.\\

If $S_8$ is taken over $d\neq1,h\neq1,g\neq1,f=1$, then
\begin{align*}
|S_8|&=\left | \underset{\alpha\in\mathbb{F}^*_{q^n}}{\sum}\underset{1\neq d,h|q^n-1}{\sum}\underset{1\neq g|x^n-1}{\sum}\frac{\mu(d)\mu(h)\mu^\prime(g)}{\phi(d)\phi(h)\Phi(g)}\underset{\chi_d,\chi_h}{\sum}\underset{\psi_g}{\sum}\chi_d(\alpha)\chi_h(\alpha^2+\alpha+1)\psi_g(\alpha)\right |\\
&\leq \underset{\underset{d,h\, square\,free}{1\neq d,h|q^n-1}}{\sum}\;\underset{\underset{g\, square\,free}{1\neq g|x^n-1}}{\sum}\frac{1}{\phi(d)\phi(h)\Phi(g)}\underset{\chi_d,\chi_h}{\sum}\underset{\psi_g}{\sum}\left |\underset{\alpha\in\mathbb{F}^*_{q^n}}{\sum}\chi_d(\alpha)\chi_h(\alpha^2+\alpha+1)\psi_g(\alpha)\right |\\
&\leq \underset{\underset{d,h\, square\,free}{1\neq d,h|q^n-1}}{\sum}\;\underset{\underset{g\, square\,free}{1\neq g|x^n-1}}{\sum}\frac{1}{\phi(d)\phi(h)\Phi(g)}\underset{\chi_d,\chi_h}{\sum}\underset{\psi_g}{\sum}\left |\underset{\alpha\in\mathbb{F}_{q^n}}{\sum}\chi_d(\alpha)\chi_h(\alpha^2+\alpha+1)\psi_g(\alpha)\right |
\end{align*}
Using theorem 3.6, we have $ \left |\underset{\alpha\in\mathbb{F}_{q^n}}{\sum}\chi_d(\alpha)\chi_h(\alpha^2+\alpha+1)\psi_g(\alpha)\right| \leq 3q^{n/2}$ and hence \\$|S_8|\leq 3q^{n/2}(2^\omega-1)^2(2^\Omega-1)$.
\vspace{.5cm}
\\ For the following, we consider 
\begin{align*}
\psi_g& = \psi_{\delta_g},\; \delta_g\in\mathbb{F}^*_{q^n}\\
\psi_f& = \psi_{\gamma_f},\; \gamma_f\in\mathbb{F}^*_{q^n}\\
\mbox{so that} \; \psi_g(\beta)&=\psi_{\delta_g}(\beta)=\psi_0(\delta_g\beta)\\
\psi_f(\beta)&=\psi_{\gamma_f}(\beta)=\psi_0(\gamma_f\beta)
\end{align*}
for $\beta\in\mathbb{F}_{q^n}\, \mbox{and}\, \psi_0$ is canonical additive character of $\mathbb{F}_{q^n}$.\\
If $S_9$ is taken over $d=h=1=g, f\neq 1$, then
\begin{align*}
|S_9| &=\left | \underset{\alpha\in\mathbb{F}^*_{q^n}}{\sum}\underset{1\neq f|x^n-1}{\sum}\frac{\mu^\prime(f)}{\Phi(f)}\underset{\psi_f}{\sum}\psi_f(\alpha^2+\alpha+1)\right|\\
&\leq \underset{\underset{f\, squarefree}{1\neq f|x^n-1}}{\sum}\frac{1}{\Phi(f)}\underset{\psi_f}{\sum}\left| \underset{\alpha\in \mathbb{F}_{q^n}}{\sum}\psi_f(\alpha^2+\alpha+1)-\psi_f(1)\right|\\
&\leq \underset{\underset{f\, squarefree}{1\neq f|x^n-1}}{\sum}\frac{1}{\Phi(f)}\underset{\psi_f}{\sum}\left\{\left| \underset{\alpha\in \mathbb{F}_{q^n}}{\sum}\psi_f(\alpha^2+\alpha+1)\right|+\left| \psi_f(1)\right|\right\}
\end{align*}
By theorem 3.4, we have$\left| \underset{\alpha\in \mathbb{F}_{q^n}}{\sum}\psi_f(\alpha^2+\alpha+1)\right|\leq q^{n/2}$ and $ |\psi_f(1)|=1$, hence \\ $ |S_9|\leq (q^{n/2}+1)(2^\Omega-1)$.\\
If $S_{10}$ is taken over $d\neq1, h=1,g=1, f\neq 1$, then 
\begin{align*}
|S_{10}|&=\left | \underset{\alpha\in\mathbb{F}^*_{q^n}}{\sum}\underset{1\neq d|q^n-1}{\sum}\underset{1\neq f|x^n-1}{\sum}\frac{\mu(d)\mu^\prime(f)}{\phi(d)\Phi(f)}\underset{\chi_d}{\sum}\underset{\psi_f}{\sum}\chi_d(\alpha)\psi_f(\alpha^2+\alpha+1)\right |\\
&\leq \underset{\underset{d\, square\,free}{1\neq d|q^n-1}}{\sum}\;\underset{\underset{f\, square\,free}{1\neq f|x^n-1}}{\sum}\frac{1}{\phi(d)\Phi(f)}\underset{\chi_d}{\sum}\underset{\psi_f}{\sum}\left |\underset{\alpha\in\mathbb{F}^*_{q^n}}{\sum}\chi_d(\alpha)\psi_f(\alpha^2+\alpha+1)\right |
\end{align*}
Using theorem 3.2, we have $\left |\underset{\alpha\in\mathbb{F}^*_{q^n}}{\sum}\chi_d(\alpha)\psi_f(\alpha^2+\alpha+1)\right |\leq 2q^{n/2}$ and hence \\ $ |S_{10}|\leq 2q^{n/2}(2^\omega-1)(2^\Omega-1)$.\\
If $S_{11}$ is taken over $d=1,h\neq1,g=1,f\neq1$, then
\begin{align*}
|S_{11}|&=\left | \underset{\alpha\in\mathbb{F}^*_{q^n}}{\sum}\underset{1\neq h|q^n-1}{\sum}\underset{1\neq f|x^n-1}{\sum}\frac{\mu(h)\mu^\prime(f)}{\phi(h)\Phi(f)}\underset{\chi_h}{\sum}\underset{\psi_f}{\sum}\chi_d(\alpha^2+\alpha+1)\psi_f(\alpha^2+\alpha+1)\right |\\
&\leq \underset{\underset{h\, square\,free}{1\neq h|q^n-1}}{\sum}\;\underset{\underset{f\, square\,lfree}{1\neq f|x^n-1}}{\sum}\frac{1}{\phi(h)\Phi(f)}\underset{\chi_h}{\sum}\underset{\psi_f}{\sum}\left |\underset{\alpha\in\mathbb{F}^*_{q^n}}{\sum}\chi_h(\alpha^2+\alpha+1)\psi_f(\alpha^2+\alpha+1)\right |\\
&\leq \underset{\underset{h\, square\,free}{1\neq h|q^n-1}}{\sum}\;\underset{\underset{f\, square\,free}{1\neq f|x^n-1}}{\sum}\frac{1}{\phi(h)\Phi(f)}\underset{\chi_h}{\sum}\underset{\psi_f}{\sum}\left |\underset{\alpha\in\mathbb{F}_{q^n}}{\sum}\chi_h(\alpha^2+\alpha+1)\psi_f(\alpha^2+\alpha+1)-\chi_h(1)\psi_f(1)\right |
\end{align*}
By applying theorem 3.6, we have  $ \left |\underset{\alpha\in\mathbb{F}^*_{q^n}}{\sum}\chi_h(\alpha^2+\alpha+1)\psi_f(\alpha^2+\alpha+1)\right |\leq(3q^{n/2}+1)$ and hence \\ $|S_{11}|\leq(3q^{n/2}+1)(2^\omega-1)(2^\Omega-1)$.\\

If $S_{12}$ is taken over $d\neq1,h\neq1,g=1,f\neq1$, then
\begin{align*}
|S_{12}|&=\left | \underset{\alpha\in\mathbb{F}^*_{q^n}}{\sum}\underset{1\neq d,h|q^n-1}{\sum}\underset{1\neq f|x^n-1}{\sum}\frac{\mu(d)\mu(h)\mu^\prime(f)}{\phi(d)\phi(h)\Phi(f)}\underset{\chi_d,\chi_h}{\sum}\underset{\psi_f}{\sum}\chi_d(\alpha)\chi_h(\alpha^2+\alpha+1)\psi_f(\alpha^2+\alpha+1)\right |\\
&\leq \underset{\underset{d,h\, square\,free}{1\neq d,h|q^n-1}}{\sum}\;\underset{\underset{f\, square\,free}{1\neq f|x^n-1}}{\sum}\frac{1}{\phi(d)\phi(h)\Phi(f)}\underset{\chi_d,\chi_h}{\sum}\underset{\psi_f}{\sum}\left |\underset{\alpha\in\mathbb{F}^*_{q^n}}{\sum}\chi_d(\alpha)\chi_h(\alpha^2+\alpha+1)\psi_f(\alpha^2+\alpha+1)\right |\\
&\leq \underset{\underset{d,h\, square\,free}{1\neq d,h|q^n-1}}{\sum}\;\underset{\underset{f\, square\,free}{1\neq f|x^n-1}}{\sum}\frac{1}{\phi(d)\phi(h)\Phi(f)}\underset{\chi_d,\chi_h}{\sum}\underset{\psi_f}{\sum}\left |\underset{\alpha\in\mathbb{F}_{q^n}}{\sum}\chi_d(\alpha)\chi_h(\alpha^2+\alpha+1)\psi_f(\alpha^2+\alpha+1)\right |
\end{align*}
By applying theorem 3.6, we have $ \left |\underset{\alpha\in\mathbb{F}_{q^n}}{\sum}\chi_d(\alpha)\chi_h(\alpha^2+\alpha+1)\psi_f(\alpha^2+\alpha+1)\right| \leq 4q^{n/2}$ and hence $|S_{12}|\leq 4q^{n/2}(2^\omega-1)^2(2^\Omega-1)$.\\

If $S_{13}$ is taken over $d=1,h=1,g\neq1,f\neq1$, then
\begin{align*}
|S_{13}|&=\left | \underset{\alpha\in\mathbb{F}^*_{q^n}}{\sum}\underset{1\neq g,f|x^n-1}{\sum}\frac{\mu^\prime(g)\mu^\prime(f)}{\Phi(g)\Phi(f)}\underset{\psi_g,\psi_f}{\sum}\psi_g(\alpha)\psi_f(\alpha^2+\alpha+1)\right |\\
&\leq \underset{\underset{g,f\, squarefree}{1\neq g,f|x^n-1}}{\sum}\frac{1}{\Phi(g)\Phi(f)}\underset{\psi_g,\psi_f}{\sum}\left |\underset{\alpha\in\mathbb{F}^*_{q^n}}{\sum}\psi_g(\alpha)\psi_f(\alpha^2+\alpha+1)\right |\\
&=\underset{\underset{g,f\, squarefree}{1\neq g,f|x^n-1}}{\sum}\frac{1}{\Phi(g)\Phi(f)}\underset{\psi_g,\psi_f}{\sum}\left |\underset{\alpha\in\mathbb{F}_{q^n}}{\sum}\psi_g(\alpha)\psi_f(\alpha^2+\alpha+1)-\psi_g(0)\psi_f(1)\right |\\
&=\underset{\underset{g,f\, squarefree}{1\neq g,f|x^n-1}}{\sum}\frac{1}{\Phi(g)\Phi(f)}\underset{\psi_g,\psi_f}{\sum}\left |\left(\underset{\alpha\in\mathbb{F}_{q^n}}{\sum}\psi_g(\alpha)\right)\left(\underset{\alpha\in\mathbb{F}_{q^n}}{\sum}\psi_f(\alpha^2+\alpha+1)\right)-\psi_g(0)\psi_f(1)\right |\\
\end{align*}
By theorem 3.1 and $|\psi_f(1)|=1$, we have $|S_{13}|\leq (2^\Omega-1)^2$.\\

If $S_{14}$ is taken over $d\neq1, h=1, f\neq1, g\neq1$, then\\
\begin{align*}
|S_{14}|&=\left | \underset{\alpha\in\mathbb{F}^*_{q^n}}{\sum}\underset{1\neq d|q^n-1}{\sum}\underset{1\neq f,g|x^n-1}{\sum}\frac{\mu(d)\mu^\prime(g)\mu^\prime(f)}{\phi(d)\Phi(g)\Phi(f)}\underset{\chi_d}{\sum}\underset{\psi_g,\psi_f}{\sum}\chi_d(\alpha)\psi_g(\alpha)\psi_f(\alpha^2+\alpha+1)\right |\\
&\leq \underset{\underset{d\, square\,free}{1\neq d|q^n-1}}{\sum}\;\underset{\underset{g,f\, square\,free}{1\neq g,f|x^n-1}}{\sum}\frac{1}{\phi(d)\Phi(g)\Phi(f)}\underset{\chi_d}{\sum}\underset{\psi_g,\psi_f}{\sum}\left |\underset{\alpha\in\mathbb{F}^*_{q^n}}{\sum}\chi_d(\alpha)\psi_g(\alpha)\psi_f(\alpha^2+\alpha+1)\right |\\
&\leq \underset{\underset{d\, square\,free}{1\neq d|q^n-1}}{\sum}\;\underset{\underset{g,f\, square\,free}{1\neq g,f|x^n-1}}{\sum}\frac{1}{\phi(d)\Phi(g)\Phi(f)}\underset{\chi_d}{\sum}\underset{\psi_g,\psi_f}{\sum}\left |\underset{\alpha\in\mathbb{F}_{q^n}}{\sum}\chi_d(\alpha)\psi_g(\alpha)\psi_f(\alpha^2+\alpha+1)\right |\\
&\leq \underset{\underset{d\, square\,free}{1\neq d|q^n-1}}{\sum}\;\underset{\underset{g,f\, square\,free}{1\neq g,f|x^n-1}}{\sum}\frac{1}{\phi(d)\Phi(g)\Phi(f)}\underset{\chi_d}{\sum}\underset{\psi_g,\psi_f}{\sum}\left |\underset{\alpha\in\mathbb{F}_{q^n}}{\sum}\chi_d(\alpha)\psi_{\delta_g}(\alpha)\psi_{\gamma_f}(\alpha^2+\alpha+1)\right |\\
&\leq \underset{\underset{d\, square\,free}{1\neq d|q^n-1}}{\sum}\;\underset{\underset{g,f\, square\,free}{1\neq g,f|x^n-1}}{\sum}\frac{1}{\phi(d)\Phi(g)\Phi(f)}\underset{\chi_d}{\sum}\underset{\psi_g,\psi_f}{\sum}\left |\underset{\alpha\in\mathbb{F}_{q^n}}{\sum}\chi_d(\alpha)\psi_0(\delta_g\alpha)\psi_0(\gamma_f(\alpha^2+\alpha+1))\right |\\
&\leq \underset{\underset{d\, square\,free}{1\neq d|q^n-1}}{\sum}\;\underset{\underset{g,f\, square\,free}{1\neq g,f|x^n-1}}{\sum}\frac{1}{\phi(d)\Phi(g)\Phi(f)}\underset{\chi_d}{\sum}\underset{\psi_g,\psi_f}{\sum}\left |\underset{\alpha\in\mathbb{F}_{q^n}}{\sum}\chi_d(\alpha)\psi_0(\delta_g\alpha+(\gamma_f(\alpha^2+\alpha+1))\right |
\end{align*}
Using theorem 3.6, we have $ \left |\underset{\alpha\in\mathbb{F}_{q^n}}{\sum}\chi_d(\alpha)\psi_0(\delta_g\alpha+(\gamma_f(\alpha^2+\alpha+1))\right |\leq 2q^{n/2}$ and hence \\ $|S_{14}|\leq 2q^{n/2}(2^\omega-1)(2^\Omega-1)^2$.\\

\noindent If $S_{15}$ is taken over $d=1,h\neq1,g\neq1,f\neq1$, then\\
$|S_{15}|=\left | \underset{\alpha\in\mathbb{F}^*_{q^n}}{\sum}\underset{1\neq h|q^n-1}{\sum}\underset{1\neq f,g|x^n-1}{\sum}\frac{\mu(h)\mu^\prime(g)\mu^\prime(f)}{\phi(h)\Phi(g)\Phi(f)}\underset{\chi_h}{\sum}\underset{\psi_g,\psi_f}{\sum}\chi_h(\alpha^2+\alpha+1)\psi_g(\alpha)\psi_f(\alpha^2+\alpha+1)\right |\\
\leq \underset{\underset{h\, square\,free}{1\neq h|q^n-1}}{\sum}\;\underset{\underset{g,f\, square\,free}{1\neq g,f|x^n-1}}{\sum}\frac{1}{\phi(h)\Phi(g)\Phi(f)}\underset{\chi_h}{\sum}\underset{\psi_g,\psi_f}{\sum}\left |\underset{\alpha\in\mathbb{F}_{q^n}}{\sum}\chi_h(\alpha^2+\alpha+1)\psi_g(\alpha)\psi_f(\alpha^2+\alpha+1)-\chi_h(1)\psi_g(0)\psi_f(1)\right |\\
\leq \underset{\underset{h\, square\,free}{1\neq h|q^n-1}}{\sum}\;\underset{\underset{g,f\, square\,free}{1\neq g,f|x^n-1}}{\sum}\frac{1}{\phi(h)\Phi(g)\Phi(f)}\underset{\chi_h}{\sum}\underset{\psi_g,\psi_f}{\sum}\left\{\left |\underset{\alpha\in\mathbb{F}_{q^n}}{\sum}\chi_h(\alpha^2+\alpha+1)\psi_g(\alpha)\psi_f(\alpha^2+\alpha+1)\right |+\left|\psi_f(1)\right|\right\}\\
\leq \underset{\underset{h\, square\,free}{1\neq h|q^n-1}}{\sum}\;\underset{\underset{g,f\, square\,free}{1\neq g,f|x^n-1}}{\sum}\frac{1}{\phi(h)\Phi(g)\Phi(f)}\underset{\chi_h}{\sum}\underset{\psi_g,\psi_f}{\sum}\left\{\left |\underset{\alpha\in\mathbb{F}_{q^n}}{\sum}\chi_h(\alpha^2+\alpha+1)\psi_{\delta_g}(\alpha)\psi_{\gamma_f}(\alpha^2+\alpha+1)\right |+ \; 1\right\}\\
\leq \underset{\underset{h\, square\,free}{1\neq h|q^n-1}}{\sum}\;\underset{\underset{g,f\, square\,free}{1\neq g,f|x^n-1}}{\sum}\frac{1}{\phi(h)\Phi(g)\Phi(f)}\underset{\chi_h}{\sum}\underset{\psi_g,\psi_f}{\sum}\left\{\left |\underset{\alpha\in\mathbb{F}_{q^n}}{\sum}\chi_h(\alpha^2+\alpha+1)\psi_0(\delta_g\alpha+(\gamma_f(\alpha^2+\alpha+1))\right |+\;1\right\}\\
\mbox{Using}\, \mbox{theorem 3.6, we have}\, \left |\underset{\alpha\in\mathbb{F}_{q^n}}{\sum}\chi_f(\alpha^2+\alpha+1)\psi_0(\delta_g\alpha+(\gamma_f(\alpha^2+\alpha+1))\right |\leq 3q^{n/2}\\
\mbox{and hence}\; |S_{15}|\leq (3q^{n/2}+1)(2^\omega-1)(2^\Omega-1)^2$\\

If $S_{16}$ is taken over $d\neq1, h\neq1,g\neq1,f\neq1$, then
\begin{align*}
|S_{16}|=  &\left | \underset{\alpha\in\mathbb{F}^*_{q^n}}{\sum}\underset{1\neq d,h|q^n-1}{\sum}\underset{1\neq f,g|x^n-1}{\sum}\frac{\mu(d)\mu(h)\mu^\prime(g)\mu^\prime(f)}{\phi(d)\phi(h)\Phi(g)\Phi(f)}\underset{\chi_d,\chi_h}{\sum}\underset{\psi_g,\psi_f}{\sum}\chi_d(\alpha)\chi_h(\alpha^2+\alpha+1)\psi_g(\alpha)\psi_f(\alpha^2+\alpha+1)\right|\\
\leq \underset{\underset{d,h\, square\,free}{1\neq d,h|q^n-1}}{\sum}&\;\underset{\underset{g,f\, square\,free}{1\neq g,f|x^n-1}}{\sum}\frac{1}{\phi(d)\phi(h)\Phi(g)\Phi(f)}\underset{\chi_d,\chi_h}{\sum}\underset{\psi_g,\psi_f}{\sum}\left |\underset{\alpha\in\mathbb{F}_{q^n}}{\sum}\chi_d(\alpha)\chi_h(\alpha^2+\alpha+1)\psi_g(\alpha)\psi_f(\alpha^2+\alpha+1)\right| \\
\leq \underset{\underset{d,h\, square\,free}{1\neq d,h|q^n-1}}{\sum}\;&\underset{\underset{g,f\, squarefree}{1\neq g,f|x^n-1}}{\sum}\frac{1}{\phi(d)\phi(h)\Phi(g)\Phi(f)}\underset{\chi_d,\chi_h}{\sum}\underset{\psi_g,\psi_f}{\sum}\left |\underset{\alpha\in\mathbb{F}_{q^n}}{\sum}\chi_d(\alpha)\chi_h(\alpha^2+\alpha+1)\psi_g(\alpha)\psi_f(\alpha^2+\alpha+1)\right| \\
\leq \underset{\underset{d,h\, square\,free}{1\neq d,h|q^n-1}}{\sum}\;&\underset{\underset{g,f\, squarefree}{1\neq g,f|x^n-1}}{\sum}\frac{1}{\phi(d)\phi(h)\Phi(g)\Phi(f)}\underset{\chi_d,\chi_h}{\sum}\underset{\psi_g,\psi_f}{\sum}\left |\underset{\alpha\in\mathbb{F}_{q^n}}{\sum}\chi_d(\alpha)\chi_h(\alpha^2+\alpha+1)\psi_{\delta_g}(\alpha)\psi_{\gamma_f}(\alpha^2+\alpha+1)\right |\\
\leq \underset{\underset{d,h\, square\,free}{1\neq d,h|q^n-1}}{\sum}\;&\underset{\underset{g,f\, square\,free}{1\neq g,f|x^n-1}}{\sum}\frac{1}{\phi(d)\phi(h)\Phi(g)\Phi(f)}\underset{\chi_d,\chi_h}{\sum}\underset{\psi_g,\psi_f}{\sum}\left |\underset{\alpha\in\mathbb{F}_{q^n}}{\sum}\chi_d(\alpha)\chi_h(\alpha^2+\alpha+1)\psi_0(\delta_g\alpha+(\gamma_f(\alpha^2+\alpha+1))\right |\\
\mbox{Using}\, &\mbox{theorem 3.6, we have}\, \left |\underset{\alpha\in\mathbb{F}_{q^n}}{\sum}\chi_d(\alpha)\chi_h(\alpha^2+\alpha+1)\psi_0(\delta_g\alpha+(\gamma_f(\alpha^2+\alpha+1))\right |\leq 4q^{n/2}\\
\mbox{and hence}\; |S_{16}|&\leq 4q^{n/2}(2^\omega-1)^2(2^\Omega-1)^2.
\end{align*}
Hence we have 
\begin{align*}
|N_{q^n}&(q^n-1,q^n-1,x^n-1,x^n-1)-\theta(q^n-1)^2\Theta(x^n-1)^2|
\leq \theta(q^n-1)^2\Theta(x^n-1)^2[ (q^{n/2}+1)(2^\omega-1)\\
&+(2q^{n/2}(2^\omega-1)^2)+(2^\Omega-1)+(q^{n/2}(2^\omega-1)(2^\Omega-1)+(2q^{n/2}+1)(2^\omega-1)(2^\Omega-1)\\
&+(3q^{n/2}(2^\omega-1)^2(2^\Omega-1))+(q^{n/2}+1)(2^\Omega-1)+(2q^{n/2}(2^\omega-1)(2^\Omega-1))\\&+(3q^{n/2}+1)(2^\omega-1)(2^\Omega-1)+(4q^{n/2}(2^\omega-1)^2(2^\Omega-1))+(2^\Omega-1)^2\\
&+(2q^{n/2}(2^\omega-1)(2^\Omega-1)^2)+(3q^{n/2}+1)(2^\omega-1)(2^\Omega-1)^2+(4q^{n/2}(2^\omega-1)^2(2^\Omega-1)^2)] 
\end{align*}
Our aim is to find pair $(q,n)$ such that $N_{q^n}(q^n-1,q^n-1,x^n-1,x^n-1)>0$\\
From above we have a sufficient condition for $N_{q^n}(q^n-1,q^n-1,x^n-1,x^n-1)>0$ is  
\begin{align*}
q^n-1>& (q^{n/2}+1)(2^\omega-1)+(2q^{n/2}(2^\omega-1)^2)+(2^\Omega-1)\\
&+(q^{n/2}(2^\omega-1)(2^\Omega-1))+(2q^{n/2}+1)(2^\omega-1)\\
&+(3q^{n/2}(2^\omega-1)^2(2^\Omega-1))+(q^{n/2}+1)(2^\Omega-1)+(2q^{n/2}(2^\omega-1)(2^\Omega-1))\\
&+(3q^{n/2}+1)(2^\omega-1)(2^\Omega-1)+(4q^{n/2}(2^\omega-1)^2(2^\Omega-1))+(2^\Omega-1)^2\\
&+(2q^{n/2}(2^\omega-1)(2^\Omega-1)^2)+(3q^{n/2}+1)(2^\omega-1)(2^\Omega-1)^2+(4q^{n/2}(2^\omega-1)^2(2^\Omega-1)^2) \\
\end{align*}
Which holds if $q^{n/2}>4.2^{2\omega+2\Omega}$. \hfill$\mathbf[4.1]$

Which our desired result. \hfill$\Box$\\

\begin{remark} This proof is not valid for $p=2$, as the theorem 3.4 is not applicable in this case, as $gcd(n,q)\neq1$ for $q=2^k$, where $k$ is a positive integer. This proof is not valid for $p=3$ also, as in this case $f(x)=x^2+x+1=(x-1)^2$, and $2|q^n-1$. So theorem 3.5 is not applicable. 
\par Since we are taking $q>3$, hence theorem 3.6 is applicable here, as no $g(x)$ of degree 2 can be expressed in the form $r(x)^q-r(x)$ in $\mathbb{F}_{q^n}[x]$. 
\end{remark}

\begin{corollary} Let $q=p^k$ where $p>3$ is prime and $k$ be a positive integer with $n|q-1$. For $n\geq35$, $(q,n)\in\mathfrak{M}$ if $p\geq11 \, \mbox{and}\, k\geq 7$.
\end{corollary}
\textbf{Proof:}
From [4.1], by calculation and using Lemma 3.9, Theorem 4.1 and Lemma 3.8
 we have $N_{q^n}(q^n-1,q^n-1,x^n-1,x^n-1)>0$ if  $q^{n/10}>4C(q^n-1)2^{2n}$ \hfill$\mathbf{[4.2]}$
\\ as by lemma 3.8, when $n|q-1$, we have $\Omega=n$.
\\ Now [4.2] is equivalent to 
\par \hspace{3cm} $log\,q>\frac{10\,log\,506.25}{n}+ 20\,log\,2$\hfill $\mathbf{[4.3]}$\\
Now for $n\geq35$, the condition holds for $q=p^k$,  $p\geq 11$ and $k\geq 7$.
\\ Hence for $n\geq35$, $(q,n)\in\mathfrak{M}$ if $p\geq11 \, \mbox{and}\, k\geq 7$.\hfill$\Box$

\vspace{.3cm}
\begin{corollary} Let $q=p^k$, where $p>3$ is prime and $k$ is a positive integer and $n$ is any positive integer such that $n\nmid q-1$. If $p\geq 5, k\geq 7$ and $n\geq72$, then $(q,n)\in\mathfrak{M}$.
\end{corollary}
\textbf{Proof:} In this case $\Omega\leq\frac{3}{4}n$ (by Lemma 3.8), then by Lemma 3.7 and Theorem 4.1, we have
\par \hspace{4cm} $q^{n/10}>4C(q^n-1)2^{\frac{3}{2}n}$ \hfill $\mathbf{[4.4]}$
\\which is equivalent to
 \par \hspace{4cm}$n> \frac{\log\,506.25}{\frac{1}{10}\,log\,q- \frac{3}{2}\,log\,2}$\hfill $\mathbf{[4.5]}$
 \\ The right hand side of [4.5] is a decreasing function of $q$ and it is positive when $q>32768$. If $q=5^7$, then the equation is true for all $n\geq 72$.
 \\ So, $(q,n)\in\mathfrak{M}$ for all $p\geq5,\, k\geq 7$ and $n\geq 72$.\hfill$\Box$\\
 
 \section{Significance of the result}
 From the result established by Anju and R.K.Sharma\cite{14} we see that there exists $\alpha\in \mathbb{F}_{q^n}$ such that $\alpha$ is primitive normal and $\alpha^2+\alpha+1$ is primitive when $q>181$. But in our result, we found that for the existence of $\alpha$ in $\mathbb{F}_{q^n}$ such that both $\alpha$ and $\alpha^2+\alpha+1$ are primitive normal, $q$ must be greater than 32768, which is quite larger than 181.\hfill$\Box$


\begin{thebibliography}{999}

\bibitem{1}
  E.Cortellini,
  Finite fields and cryptology,
  \emph{Computer Science \\Journal of Moldova}.
  vol.11, no.2(32),
  2003.
\bibitem{2}
   L.Fu and D.Q.Wan,
   A class of incomplte character sums,
   \emph{Q.J.Math.Soc},
   \textbf{43}, (1968) 21-39.
\bibitem{3}   
   S.D.Cohen,
   Consecutive primitive roots in a finite field,
   \emph{Proc. Amer. Math. Soc.},
   \textbf{93}(2) (1985) 189-197.
\bibitem{4}
   S.D.Cohen and S.Huczynska,
   The primitive normal basis theorem without a computer,
   \emph{J. Lond. Math. Soc.}
   \textbf{67}(1) (2003) 41-56
\bibitem{5} 
   S.D.Cohen and S.Huczynska,
   The strong primitive normal basis theorem,
   \emph{Acta. Arith.}
   \textbf{143}(4) (2010) 299-332
\bibitem{6}
   S.D.Cohen,
   Pairs of primitive elements in fields of even order,
   \emph{Finite Fields Appl.},
   \textbf{28} (2014) 22-42
\bibitem{7}
   D.Wan,
   Generators and irreducible polynomials over finite fields,
   \emph{Math. Comp.}
   \textbf{66}(219) (1997) 1195-1212
\bibitem{8}
   F.N.Castro and C.J.Moreno,
   Mixed exponential sums over finite fields,
   \emph{Proc. Amer. Math. Soc.} ,
   \textbf{128}(9) (2000) 2529-2537
\bibitem{9}
   G.Kapetankis,
   Normal bases and primitive elements over finite fields,
   \emph{Finite Fields Appl.}  
   \textbf{26}(2014) 123-143 
\bibitem{10}
   T.Garefalakis and G.Kapetanakis,
   On the existence of primitive completely normal bases of finite fields,
   \emph{J. Pure Appl. Algebra}
   (2018)  
\bibitem{11}
   H.W.Lenstra,Jr. and R.J.Schoof,
   Primitive Normal Bases for Finite Fields,
   \emph{Math. Comp.}
   \textbf{48} (1987) 217-231
\bibitem{12}
   L.Carlitz,
   Primitive roots in a finite fields,
   \emph{Trans. Amer. Math. Soc.}
   \textbf{73}(3) (1952) 314-318
\bibitem{13}
   R.Lidl and H.Niederreiter,
   \emph{Finite Fields},
   2nd edn. 
   (Cambridge \\University Press, Cambridge, 1997)
\bibitem{14}
   Anju and R.K.Sharma,
   On primitive normal elements over finite fields,
   \emph{Asian-Eur. J. Math.} ,
   \textbf{11}(2) (2018)   
\bibitem{15}
   G.James and M.Liebeck,
   \emph{Representations and Characters of Groups},
   2nd edn.
   (Cambridge University Press, Cambridge, 2001)
\bibitem{16}
   P.P.Wang,X.W.Cao and R.Q.Feng,
   On the existence of some specific elements in finite fields of characteristic 2
   \emph{Finite Fields Appl.},
   \textbf{18}(4) (2012) 800-8013   
\bibitem{17}
   Q.Liao, J.Li and K.Pu,
   On the existence for some primitive elements in finite fields,
   \emph{Chin. Ann. Math.},
   B\textbf{37} (2016) 259-266 
\bibitem{18}
   T.Tian and W.F.Qi,
   Primitive normal elements and its inverse in finite fields,
   \emph{Acta. Math. Sinica}(Chin. Ser.),
   \textbf{49}(3) (2006) 657-668                         
\end{thebibliography}
 \end{document}